\documentclass{amsart}

\usepackage{graphicx}
\usepackage{longtable}
\newtheorem{thm}{Theorem}
\newtheorem{defn}{Definition}
\newtheorem{cor}{Corollary}
\newtheorem{rem}{Remark}

\begin{document}

\title{Series of Reciprocal Powers of $k-$almost Primes}

\author{Richard J. Mathar}
\urladdr{http://www.strw.leidenuniv.nl/~mathar}
\email{mathar@strw.leidenuniv.nl}
\address{Leiden Observatory, Leiden University, P.O. Box 9513, 2300 RA Leiden, The Netherlands}


\date{\today}

\begin{abstract}
Sums over inverse $s$-th powers of semiprimes and
$k$-almost primes are transformed into sums over products
of powers of ordinary prime zeta functions. Multinomial coefficients
known from the cycle decomposition of permutation groups play the role
of expansion coefficients.
Founded on a known convergence acceleration for the ordinary
prime zeta functions, the sums and first derivatives are tabulated with high
precision for indices $k=2,\ldots,6$ and
integer powers $s=2,\ldots,8$.
\end{abstract}

\maketitle
\section{Overview}
Series over rational polynomials evaluated at integer arguments contain
sub\-series summing over integers classified by the
count of their prime factors. The core example is the Riemann zeta function $\zeta$
which accumulates the prime zeta function $P_1$ plus what we shall define
the almost-prime zeta functions $P_k$ (Section \ref{sec.2}).
The central observation of this manuscript is that the almost-prime zeta functions
are combinatorial sums over the prime zeta function (Section \ref{sec.3}).
Since earlier work by Cohen, Sebah and Gourdon has pointed at efficient
numerical algorithms to compute $P_1$, series over reciprocal almost-primes---which
may suffer from slow convergence in their defining format---may
be computed efficiently by reference to the $P_1$.

In consequence, any
converging series over the positive integers which has a Taylor expansion in reciprocal
powers of these integers splits into $k$-almost prime components.
Section \ref{sec.4} illustrates this for the most basic formats.

Number theory as such will not be advanced. The meromorphic landscape of the
prime zeta functions as a function of their main variable appears to be more
complicated than what is known for their host, the Riemann zeta function;
so only some remarks on the calculation of first derivatives are dropped.

\section{Prime Zeta Function} \label{sec.2}
\begin{defn}
The prime zeta function $P(s)$ is the sum over the reciprocal $s$-th powers
of the prime numbers $p$
\cite{FrobergBIT8,Fincharxiv06,MerrifieldPRSL33}
\begin{equation}
P(s)\equiv \sum_{p}^\infty \frac{1}{p^s}
;
\quad
\Re s>1
.
\label{eq.P1def}
\end{equation}
\end{defn}
\begin{rem}
The primes are represented by sequence
A000040
in the Online-En\-cy\-clo\-pedia of Integer Sequences
\cite{EIS},
and we will adopt the nomenclature that a letter A followed by a 6-digit number
points at a sequence in this database.
Accurate values $P(s)$ for $s\le 9$ are provided by the sequences
A085548,
A085541,
and
A085964--A085969.
\end{rem}
Table \ref{tab.P} complements these by using \cite{BachIPL62,SebahGourdon,Cohen}
\begin{equation}
P(s)=\sum_{p\le M} \frac{1}{p^s}+\sum_{n=1}^\infty
\frac{\mu(n)}{n}\log P(M,sn)
\label{eq.SG}
\end{equation}
for a suitably large prime $M$, where $\mu$ is the M\"obius
function \cite[(24.3.1)]{AS},
where $\zeta$ is the Riemann zeta function
\cite{Chengarxiv07,KolbigMathComp64_449}\cite[(9.5)]{GR},
and
\begin{equation}
P(M,s)\equiv \zeta(s)\prod_{p\le M}\left(1-p^{-s}\right)
\end{equation}
an associated definition of a partial product.

\begin{table}
\caption{The Prime Zeta Function of some integer arguments.
In a style adopted from \cite {AS},
optional trailing parentheses contain
an additional power of 10\@. Example: The number $3.45\times 10^{-3}$ may
appear as $0.00345$ or as $3.45(-3)$ or $.345(-2)$.
Trailing dots indicate that more digits are chopped off, not rounded, at
the rightmost places.
}
\begin{tabular}{|l|l|}
\hline
$s$ & $P(s)$ \\
\hline
10 &  .9936035744369802178558507001477394163018725452852033205535666(-3)
$\ldots$
\\
11 &  .4939472691046549756916217683343987121559397009604952181866074(-3)
$\ldots$
\\
12 &  .2460264700345456795266485921650809279799322679473231921741459(-3)
$\ldots$
\\
13 &  .1226983675278692799054887924314033239147428525577690135256528(-3)
$\ldots$
\\
14 &  .6124439672546447837750803429987454197282126872378013541885303(-4)
$\ldots$
\\
15 &  .3058730282327005256755462931371262800130114525389809330765981(-4)
$\ldots$
\\
16 &  .1528202621933934418080192641189055977466126987760393110788060(-4)
$\ldots$
\\
17 &  .7637139370645897250904556043939762017569839042162662520251345(-5)
$\ldots$
\\
18 &  .3817278703174996631227515316311091361624942636382614195748077(-5)
$\ldots$
\\
19 &  .1908209076926282572186179987969776618145616195068986381165765(-5)
$\ldots$
\\
20 &  .9539611241036233263528834939770057955700555885822134364992986(-6)
$\ldots$
\\
21 &  .4769327593684272505083726618818876106041908102543778311286107(-6)
$\ldots$
\\
22 &  .2384504458767019281263116852015955086787325069914706736138961(-6)
$\ldots$
\\
23 &  .1192199117531882856160246453383398577108304116591413496750467(-6)
$\ldots$
\\
24 &  .5960818549833453297113066655008620131146582480117715598724992(-7)
$\ldots$
\\
25 &  .2980350262643865876662659401778145949592778827831139923166297(-7)
$\ldots$
\\
26 &  .1490155460631457054345907739442373384026574094003717094826319(-7)
$\ldots$
\\
27 &  .7450711734323300780164546124093693349559346148927152517177541(-8)
$\ldots$
\\
28 &  .3725334010910506351833912287693071753007133176180958755544001(-8)
$\ldots$
\\
29 &  .1862659720043574907522145113353601172883347161316571090677067(-8)
$\ldots$
\\
30 &  .9313274315523019206770664589654477590951135917359845054142758(-9)
$\ldots$
\\
31 &  .4656629062865372188024756168924550748371110904683071460848803(-9)
$\ldots$
\\
32 &  .2328311833134403149136721429290134383956012839353792695549588(-9)
$\ldots$
\\
33 &  .1164155017134526496600716286019717301900642951759025845555278(-9)
$\ldots$
\\
34 &  .5820772087563887361296110329279891461135544371461870778050157(-10)
$\ldots$
\\
35 &  .2910385044412396334030528313212481809718543835635609582681388(-10)
$\ldots$
\\
36 &  .1455192189083022590216132905087468529073564045445529854681582(-10)
$\ldots$
\\
37 &  .7275959835004541439158484817671131286009806802799652884873652(-11)
$\ldots$
\\
38 &  .3637979547365416297743239172591421915260411920812352913301657(-11)
$\ldots$
\\
39 &  .1818989650303757224685763903905856620025233007531508321634193(-11)
$\ldots$
\\
\hline
\end{tabular}
\label{tab.P}
\end{table}

Moments are listed
in Table \ref{tab.Pmom}. The first line,
\begin{equation}
\sum_{s=1}^\infty \frac{1}{s}P(s)
 = -\sum_p \left[\log\left(1-\frac{1}{p}\right)+\frac{1}{p}\right]
 = \lim_{s\to 1}[ \zeta(s)-\sum_p \frac{1}{p^s}]
\label{eq.L11}
,
\end{equation}
is A143524 \cite{EIS,Cohen}.
\begin{table}
\caption{
Moments of the Prime Zeta function.
The difference between the
first two lines has been quoted by Chan \cite{ChanAA123}.
}
\begin{tabular}{|l|l|}
\hline
$u$ & $\sum_{s=2}^\infty P(s)/s^u = \sum_{s=2}^\infty \sum_p 1/(s^u p^s)$ \\
\hline
1 & .31571845205389007685108525147370657199059268767872439261370
$\ldots$
\\
2 & .139470639611308681803077672937394594285963216485548292900601
$\ldots$
\\
3 & .646081378884568610840336242581546273479300400893390674692107(-1)
$\ldots$
\\
4 & .307990601127799353645768912808351277991590439776491610562576(-1)
$\ldots$
\\
5 & .149414104078444272503589211881003770331602105485328925740449(-1)
$\ldots$
\\
6 & .732763762441199457647551708167392749211019881968761942480213(-2)
$\ldots$
\\
\hline
\end{tabular}
\label{tab.Pmom}
\end{table}
The total of all values in Table \ref{tab.Pmom} is
\begin{eqnarray}
\sum_{u=1}^\infty \sum_{s=2}^\infty \frac{1}{s^u}P(s)
&=&
\sum_{s=2}^\infty \sum_p \frac{1}{(s-1)p^s}=
-\sum_p \frac{ \log(1-1/p)}{p}
\nonumber
\\
&=& 0.58005849381391172358283349737677118691587319037\ldots
\label{eq.pmoms}
\end{eqnarray}

The first derivatives $dP(s)/ds$ of (\ref{eq.P1def}) are
evaluated as the first derivatives of (\ref{eq.SG}) \cite[(3.3.6)]{AS},
\begin{eqnarray}
P'(s)
&=&
-\sum_{p} \frac{\log p}{p^s}
=
-\sum_{p\le M} \frac{\log p}{p^s}+\sum_{n=1}^\infty
\mu(n) \frac{P'(M,sn)}{P(M,sn)};
\label{eq.Pprime}
\\
\frac{P'(M,s)}{P(M,s)}
&=&
\frac{\zeta'(s)}{\zeta(s)}
+
\sum_{p\le M}\frac{\log p}{p^s(1-p^{-s})}.
\label{eq.PMprime}
\end{eqnarray}
Here, primes denote derivatives with respect to the main argument, which is the
second argument for the case of $P(.,.)$.
Table \ref{tab.Pp} shows some of them for selected small integer $s$.
\begin{rem}
The relation $-\zeta'(s)/\zeta(s) = \sum_p \log(p)/(p^s-1)$ is a gateway
to acceleration of series involving logarithmic numerators. Examples
are  growth rate coefficients
of unitary square-free divisors \cite{SuryanJRAM276},
\begin{gather}
\sum_p \frac{(2p+1)\log p}{(p+1)(p^2+p-1)}
=
\sum_p \log p\Big[\frac{2}{p^2-1}
-\frac{3}{p^3-1}
+\frac{4}{p^4-1}
-\frac{10}{p^5-1}
+\frac{18}{p^6-1}
\nonumber
\\
\quad
-\frac{28}{p^7-1}
+\frac{40}{p^8-1}
-\frac{72}{p^9-1}
+\cdots
\Big]
= 0.748372333429674\ldots
\end{gather}
or of unitary cube-free divisors,
\begin{gather}
\sum_p \frac{(4p^2-2-p)\log p}{(p^2-1)(p^2+p-1)}
=
\sum_p \log p\Big[\frac{4}{p^2-1}
-\frac{5}{p^3-1}
+\frac{7}{p^4-1}
-\frac{17}{p^5-1}
+\frac{31}{p^6-1}
\nonumber
\\
\quad
-\frac{48}{p^7-1}
+\frac{69}{p^8-1}
-\frac{124}{p^9-1}
+\cdots
\Big]
= 1.647948081159756\ldots
\end{gather}

\end{rem}

\begin{table}
\caption{The first derivative of the Prime Zeta Function at some integer arguments $s$.
The value in the first line differs from Cohen's value \cite{Cohen} after 42 digits.
}
\begin{tabular}{|l|l|}
\hline
$s$ & $P'(s)$ \\
\hline
2 & -4.930911093687644621978262050564912580555881263464682907133271(-1)
$\ldots$
\\
3 & -1.507575555439504221798365163653429195755011615306893318187976(-1)
$\ldots$
\\
4 & -6.060763335077006339223098370971337840638287746125984399112768(-2)
$\ldots$
\\
5 & -2.683860127679835742218751329245015994333014955355822812481980(-2)
$\ldots$
\\
6 & -1.245908072279999152702779277468997004091135047157587587410933(-2)
$\ldots$
\\
7 & -5.940689039148196142550592829016609019368189505929351075166813(-3)
$\ldots$
\\
8 & -2.879524708729247391346028423857334064998983761675865841067618(-3)
$\ldots$
\\
9 & -1.410491921424531291554196456308199977901657131693496192836500(-3)
$\ldots$
\\
10 & -6.956784473446204802000701977708415913844863703329838954712256(-4)
$\ldots$
\\
11 & -3.446864256305149016520798301347221055148509398720732052598028(-4)
$\ldots$
\\
12 & -1.712993524462175657532493112138275372004981118241302276420951(-4)
$\ldots$
\\
13 & -8.530310916711056635208876017215691972617326615054214472499073(-5)
$\ldots$
\\
14 & -4.253630557412291035554757415368617516720893534438843631304558(-5)
$\ldots$
\\
15 & -2.122979056274934599669348621302375720453112762226994727150844(-5)
$\ldots$
\\
16 & -1.060211861676127903320578231686279299852887328732516230264968(-5)
$\ldots$
\\
17 & -5.296802557643848074496697331902062291354582070044729659167083(-6)
$\ldots$
\\
18 & -2.646982787802997352263261854182101806956865359404392741570106(-6)
$\ldots$
\\
19 & -1.323018648512292735443206851957658773372595611301942028763990(-6)
$\ldots$
\\
20 & -6.613517594172600210891457029052100560435779681754585968634604(-7)
$\ldots$
\\
21 & -3.306233614825208657730023089591286331399889373034673500892396(-7)
$\ldots$
\\
22 & -1.652941753425972669328543237067224505237754606957895309294978(-7)
$\ldots$
\\
23 & -8.264125267365738127779160862943622945349018242909935277495509(-8)
$\ldots$
\\
24 & -4.131868136465068742054546598016395808846430264964223880702300(-8)
$\ldots$
\\
25 & -2.065869236367122379085627896761801003317031762395199472960982(-8)
$\ldots$
\\
26 & -1.032913007669833840610060139473968867069681131430230354135307(-8)
$\ldots$
\\
27 & -5.164493003519525183949097124602884025560927707601929381567114(-9)
$\ldots$
\\
28 & -2.582222490193098373680412778703362364263231350157972835761708(-9)
$\ldots$
\\
29 & -1.291103241249637065884459649285079993129747526229207669548247(-9)
$\ldots$
\\
\hline
\end{tabular}
\label{tab.Pp}
\end{table}

\section{Zeta Functions of Almost-primes} \label{sec.3}
\subsection{Nomenclature}

We generalize the notation, and define the $k$-almost prime zeta functions
by summation over inverse powers of $k$-almost primes $q_j$,
\begin{defn}
\begin{equation}
P_k(s)\equiv \sum_{j=1}^\infty \frac{1}{q_j^s}.
\label{eq.P1}
\end{equation}
In slight violation of the \emph{almost}-terminology, the Prime Zeta Function is incorporated
as just one special case,
\begin{equation}
P_1(s)\equiv P(s).
\end{equation}
\end{defn}
\begin{rem}
The sequence $q_j$ is given
by the primes
{A000040} if $k=1$,
by the semiprimes
{A001358} if $k=2$,
by the 3-almost primes
{A014612} if $k=3$,
by the 4-almost primes
{A014613} if $k=4$,
by the 5-almost primes
{A014614} if $k=5$,
by the 6-almost primes
{A046306} if $k=6$,
by the 7-almost primes
{A046308} if $k=7$
etc.
$P_2(2)$ is
A117543,
and $P_3(2)$ is
A131653.
\end{rem}
\begin{rem}
One step further defines prime multi-zeta functions of the form
\begin{equation}
P(s_1,s_2,\ldots,s_k)\equiv \sum_{p_1,p_2,\ldots ,p_k=2,3,5\ldots}\frac{1}{p_1^{s_1}p_2^{s_2}\cdots p_k^{s_k}}.
\end{equation}
If the exponents are restricted to a common number $s$, they retrieve the information
of the $P_k(s)$ in slightly entangled form, for example $P(s,s)=2P_2(s)-P(2s)$ \cite{BradleyDevMath14}.
\end{rem}
Each
integer $n>0$ is either a member of the set $\{1\}$,
or of the set of primes, or of the set of semi-primes, etc. These
disjoint sets are labeled by the sum of the exponents of the prime
number factorizations of their members,
\begin{equation}
\Omega(n)=k;\quad \Omega(1)\equiv 0.
\end{equation}
The Riemann zeta function
may be partitioned into sums over the almost-prime zeta functions,
\begin{equation}
\zeta(s)
=
\sum_{n=1}^\infty \frac{1}{n^s}
=
1+\sum_{k=1}^\infty \sum_{\substack{n=1\\ \Omega(n)=k}}^\infty \frac{1}{n^s}
=
1+\sum_{k=1}^\infty P_k(s).
\label{eq.zetpart}
\end{equation}
\begin{rem}
$\zeta(s)$ may be taken from \cite[(Table 23.3)]{AS} for $s\le 42$,
or from 
A013661,
A002117,
A013662
--
A013678
while $s\le 20$, or from the link to ``Recent additions of tables'' in Plouffe's database \cite{Plouffe}
while $s\le 99$.
$\Omega(n)$ is tabulated in A001222.
\end{rem}

\subsection{Numerical Scheme}

Table \ref{tab.Pk} of $P_k(s)$ is deduced from
\begin{thm}\label{thm.main}
\cite{KnopfmacherCN173}
\begin{gather}
P_2(s)=\frac{P(2s)+P^2(s)}{2!};
\label{eq.P2}
\\
P_3(s)=\frac{2P(3s)+3P(2s)P(s)+P^3(s)}{3!};
\\
P_4(s)=\frac{6P(4s)+8P(3s)P(s)+3P^2(2s)+6P(2s)P^2(s)+P^4(s)}{4!};
\\
P_k(s)=\frac{1}{k!}\sum_{\substack{k_1+2k_2+3k_3+\dotsb+kk_k=k\\ k_k\ge 0}}(k;k_1k_2\ldots k_k)^*
P^{k_1}(s)P^{k_2}(2s)\dotsm P^{k_k}(ks),
\label{eq.Pks}
\end{gather}
utilizing values of the prime zeta function as discussed above
and summing over the partitions $\pi(k)$ of $k$ with
weight coefficients
\begin{equation}
(k;k_1k_2\ldots k_k)^*=k!/\prod_{m=1}^k(m^{k_m}k_m!)
\label{eq.m2def}
\end{equation}
of Table 24.2 in \cite{AS}
(multinomials $M_2$) and A036039 or A102189.
\end{thm}

\begin{proof}
(\ref{eq.Pks}) is the main result of the paper. For small $k$, explicit
verification can be done along Price's 
\cite{PriceAMM53}
construction, where the $k$-almost primes fill triangular,
($k=2$), tetrahedral ($k=3$) etc.\ sections of a
$k$-dimensional Euclidean lattice labeled by the prime numbers
along its Cartesian axes \cite{BennettAlgU32}. The case $P_1(s)=P(s)$ just repeats
the definition (\ref{eq.P1}). The case (\ref{eq.P2}) accumulates
in $P(2s)$ the sum over the squares, and in $P^2(s)$---with the binomial expansion---
again the sum over the squares and twice the sum over products of distinct primes.
After division through $2!$, each semiprime is effectively represented once.

The generic proof follows through induction:
the terms of the right hand side of (\ref{eq.Pks})
contain the factor
$(k;k_1k_2,\ldots k_k)^*$, which is the number of distinct permutations
with $k_m$ cycles of length $m$ for $m=1,2\ldots k$
\cite[(p.\ 123)]{Berge}\cite[(24.1.2)]{AS}.
The right hand side is the cycle index
\begin{equation}
Z(S_k)\equiv \frac{1}{k!}\sum_{\substack{k_1+2k_2+3k_3+\dotsb+kk_k=k\\ k_k\ge 0}}(k;k_1k_2\ldots k_k)^*
P^{k_1}(s)P^{k_2}(2s)\dotsm P^{k_k}(ks)
\end{equation}
of the symmetric group $S_k$
\cite[(2.2.5)]{Harary} with $P(ms)$ substituted for the indeterminates
of cycle length $m$. Skipping any interpretation within a Redfield-P\'olya
symmetry, its recurrence
\begin{equation}
Z(S_k)=\frac{1}{k}\sum_{j=1}^{k} P(js) Z(S_{k-j})
\label{eq.cyclRec}
\end{equation}
is already established \cite[(2.2.9)]{Harary}.
This matches precisely the recurrence on the left hand side which generates
$P_k$ by a combination of products of lower-indexed almost-primes,
\begin{equation}
P_k(s)=\frac{1}{k}\sum_{j=1}^{k} P(js) P_{k-j}(s).
\end{equation}
This recurrence is valid because each $k$-almost prime which appears on
the left hand side of this equation can be generated in $k$ ways by
a product of the form $P(js)P_{k-j}(s)$: in $\omega(k)$ ways by splitting
off a prime number and multiplication with a number of the sum in $P_{k-1}(s)$,
for each divisor of the $k$-almost prime which is a square of some prime in
addition by multiplication of the square with a term in the sum in $P_{k-2}(s)$, 
and so on for divisors that are cubes of some prime etc\@.
\end{proof}

\begin{figure}
\includegraphics[scale=0.7]{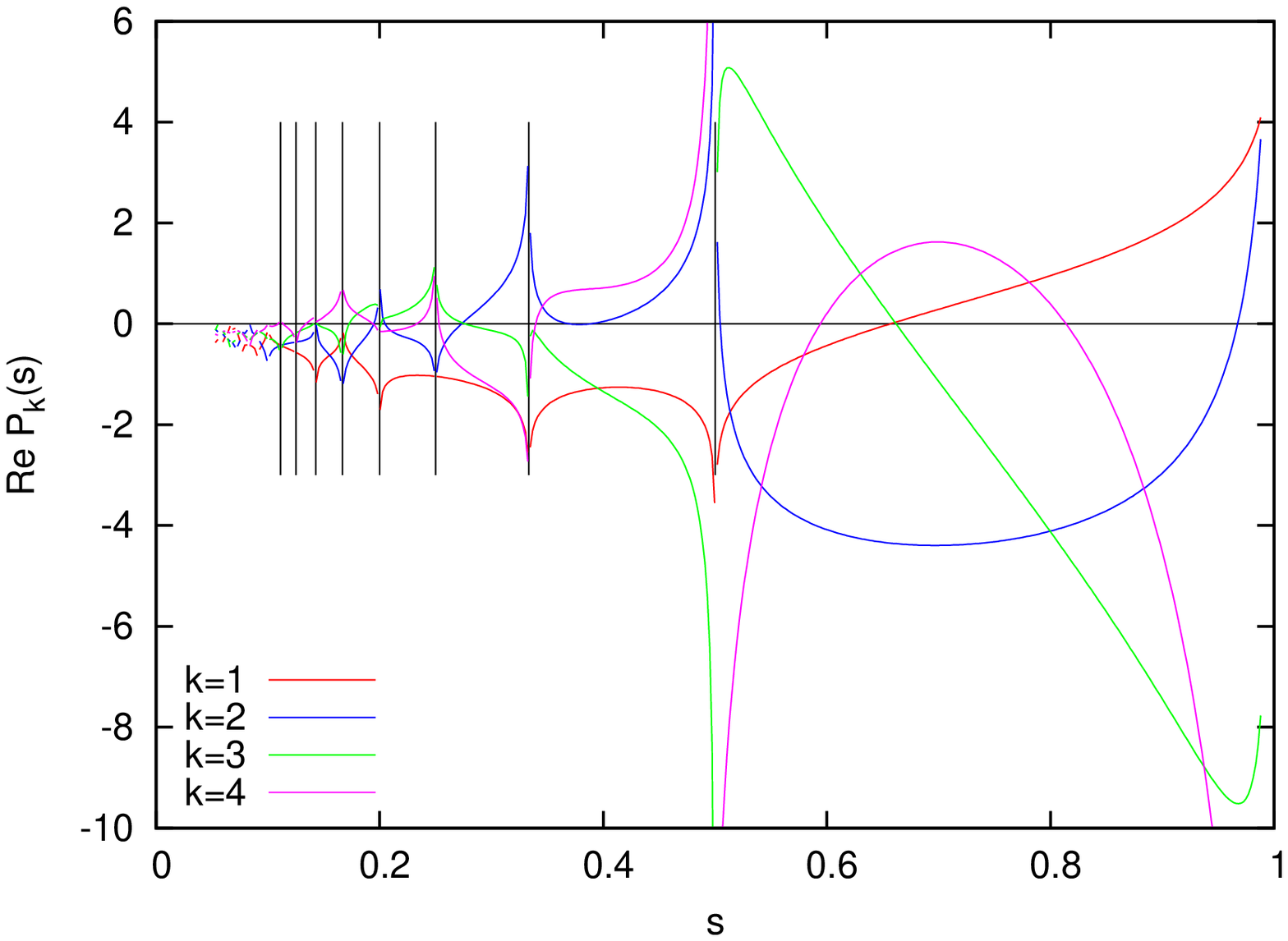}
\caption{Structure of poles and logarithmic singularities:
$\Re P_k(s)$ on the real line between 0 and 1.
}
\label{fig.frob}
\end{figure}

\begin{table}
\caption{Almost-prime zeta functions $P_k(s)$ at small integer arguments $s$.
}
\begin{tabular}{|r|r|l|}
\hline
$k$ & $s$ & $P_k(s)$ \\
\hline
2 & 2 & 1.407604343490233882227509254138772537749192760048802639241489(-1)
$\ldots$
\\
2 & 3 & 2.380603347277195967869595585283620062893217848034845684562765(-2)
$\ldots$
\\
2 & 4 & 4.994674468637339635276874049579289322502057848230867728509096(-3)
$\ldots$
\\
2 & 5 & 1.136012424856354766515556190735772665693748056026108556151424(-3)
$\ldots$
\\
2 & 6 & 2.687071675614096324217387396140875535798787447125719642936101(-4)
$\ldots$
\\
2 & 7 & 6.493314175691145578854061507836714519989975167152833237216508(-5)
$\ldots$
\\
2 & 8 & 1.588851988525958888572372095351879234527858971327233300748108(-5)
$\ldots$
\\
\hline
3 & 2 & 3.851619298269464091283792262806039543890016747838157193719155(-2)
$\ldots$
\\
3 & 3 & 3.049362082334312946748098847079302999848694548619577993637287(-3)
$\ldots$
\\
3 & 4 & 3.144274968329417421821246641907192073071706953574340102524412(-4)
$\ldots$
\\
3 & 5 & 3.557725337068269111888017622799305930206716602282958084700573(-5)
$\ldots$
\\
3 & 6 & 4.201275533960671214387834295923202794959720951879928447823823(-6)
$\ldots$
\\
3 & 7 & 5.073887994515979227127878654920650441797124899213497891489656(-7)
$\ldots$
\\
3 & 8 & 6.206813624161469945551964458392656691354524774013736254471908(-8)
$\ldots$
\\
\hline
4 & 2 & 1.000943620148325082041084351808525466652473851036634849174401(-2)
$\ldots$
\\
4 & 3 & 3.839045346157269074628008425162843300890790106333110559279434(-4)
$\ldots$
\\
4 & 4 & 1.967963362818191467940855961573410955099879950428233958199339(-5)
$\ldots$
\\
4 & 5 & 1.112105498394147042065416843932409614810339288829649464410277(-6)
$\ldots$
\\
4 & 6 & 6.564866966272364593992630942917565336419279606893244510336952(-8)
$\ldots$
\\
4 & 7 & 3.964020093813893558567748245375642870705531500802782383854435(-9)
$\ldots$
\\
4 & 8 & 2.424542067198129719213221460827573885198415530509018971808441(-10)
$\ldots$
\\
\hline
5 & 2 & 2.545076168069302058221776985605516223099431333404435645812102(-3)
$\ldots$
\\
5 & 3 & 4.808940110832567973019045453666287670709263774628437825310148(-5)
$\ldots$
\\
5 & 4 & 1.230321747728495443208363890849979176133316153001252000782537(-6)
$\ldots$
\\
5 & 5 & 3.475459860092756789327837058184938607371038782365590655202548(-8)
$\ldots$
\\
5 & 6 & 1.025765593034930528602801778254441805529589443170682127946500(-9)
$\ldots$
\\
5 & 7 & 3.096892760390829520074774635913281487435638340106655991757766(-11)
$\ldots$
\\
5 & 8 & 9.470868287557099531707292414885674011014480791459591532172493(-13)
$\ldots$
\\
\hline
6 & 2 & 6.410338528642807128627067320846767912178898525413482485882042(-4)
$\ldots$
\\
6 & 3 & 6.014928780179108948186295382866155223857573977633433216567001(-6)
$\ldots$
\\
6 & 4 & 7.689936414615761724089452218863183766130186637100554867671046(-8)
$\ldots$
\\
6 & 5 & 1.086086563383684175516301346294503063069865950783005109331817(-9)
$\ldots$
\\
6 & 6 & 1.602759442759730816486790711011122437362717548888087849039669(-11)
$\ldots$
\\
6 & 7 & 2.419447563344257658856721989501609732074582033968956770368419(-13)
$\ldots$
\\
6 & 8 & 3.699557937592796079194223671580493130944594965925089150589532(-15)
$\ldots$
\\
\hline
\end{tabular}
\label{tab.Pk}
\end{table}
The simple pole of $\zeta(s)$ at $s-1$ with Stieltjes constants $\gamma_j$
\cite{Coffeyarxiv07,ApostolMathComp44,KreminskiMComp72},
\begin{equation}
\zeta(s)=\frac{1}{s-1}+\gamma+\sum_{j=1}^\infty \frac{(-1)^j}{j!} \gamma_j(s-1)^j,
\end{equation}
is associated with a logarithmic singularity of $P(s)$ at $s=1$
and singularities of $P(s)$ on the real line between $s=0$
and $s=1$ where $s$ is the inverse of a square-free integer \cite{FrobergBIT8}.
If $k>1$, the $k$-almost zeta functions inherit these and add more by the
mechanism evident from the multipliers in Theorem \ref{thm.main}.
$P_2(s)$ in (\ref{eq.P2}), for example, inherits
singularities at $1/2$, $1/3$, $1/5$,
$1/6$, $1/7$, $1/10$ etc.\ from the term $P^2(s)$, and singularities
at $1/2$, $1/4$, $1/6$, $1/10$, $1/12$ etc\@.  from the term $P(2s)$,
illustrated in Fig.\ \ref{fig.frob}.

\begin{rem}
Sums over odd indices are \cite{SebahGourdonPi,AokiMA}:
\begin{gather}
\sum_{k=1}^\infty P_{2k-1}(2)=\frac{\pi^2}{20};\qquad
\sum_{k=1}^\infty P_{2k-1}(2s)=\frac{\zeta^2(2s)-\zeta(4s)}{2\zeta(2s)}.
\end{gather}
Sums over even indices follow from there as the complement in (\ref{eq.zetpart}).
\end{rem}

\begin{figure}
\includegraphics[scale=0.7]{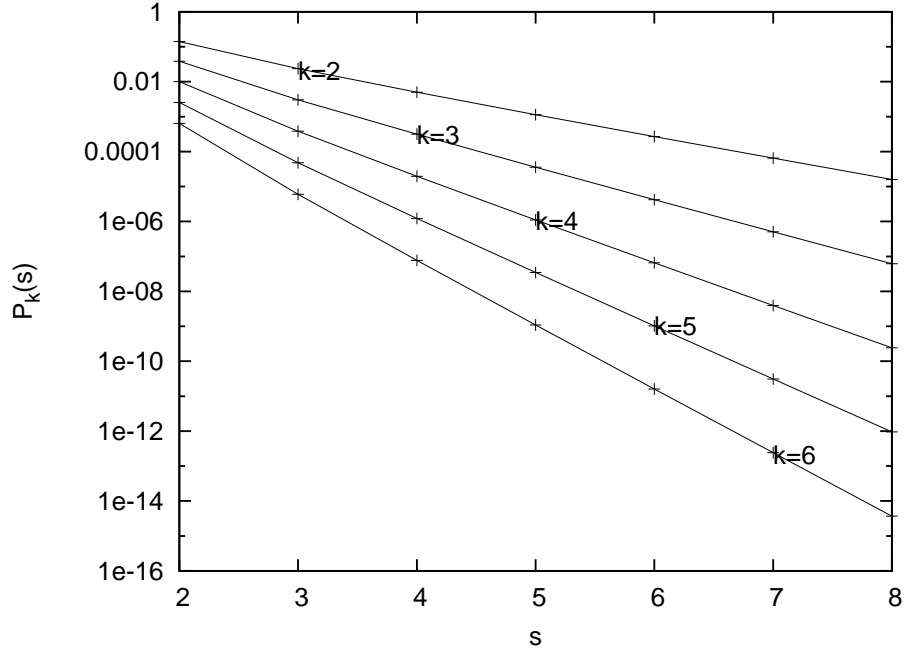}
\caption{A synopsis of table \ref{tab.Pk} on a semi-logarithmic
scale, indicating that the $P_k(s)$ fall off approximately exponentially
$P_k(s)\sim 2^{-ks}$
as $s\to \infty$ along the real $s$-axis.
}
\end{figure}

\begin{cor}
The product rule \cite[(3.3.3)]{AS} applied to (\ref{eq.Pks}) yields
the first derivative
\begin{multline}
P_k'(s)
=
\frac{1}{k!}
\sum_{\substack{k_1+2k_2+3k_3+\dotsb+kk_k=k\\ k_k\ge 0}}
(k;k_1k_2\ldots k_k)^* P^{k_1}(s)P^{k_2}(2s)\dotsm P^{k_k}(ks)
\\
\times
\left[
\frac{k_1 P'(s)}{P(s)}
+
\frac{2k_2 P'(2s)}{P(2s)}
+\dotsb
+
\frac{kk_k P'(ks)}{P(ks)}
\right]
.
\label{eq.Pksprime}
\end{multline}
\end{cor}
Numerical evaluation yields Table \ref{tab.Pkpri}.
The underivative $\int_x^\infty P_k(s)ds= \sum_i 1/[p_i^x \log p_i]$ has been evaluated
numerically for $k=1$ by Cohen \cite{Cohen,MatharArxiv0811b}.

\begin{table}
\caption{First derivatives $P_k'(s)$.
$P_1'(s)$ is Table \ref{tab.Pp}.
}
\begin{tabular}{|r|r|l|}
\hline
$k$ & $s$ & $P_k'(s)$ \\
\hline
2 & 2 & -2.836068154079806522242582225482783360793505782378140134111118(-1)
$\ldots$
\\
2 & 3 & -3.880586902399322336692460182658731189674126497852135409952156(-2)
$\ldots$
\\
2 & 4 & -7.545896694085315206907970667196350193021639416359677908501448(-3)
$\ldots$
\\
2 & 5 & -1.655293105240617640761013220868097514629396965866818280047120(-3)
$\ldots$
\\
2 & 6 & -3.839769424635045625755371800119456139211625839983225237893684(-4)
$\ldots$
\\
2 & 7 & -9.174798062469143952346462602570231884264128810420282594760147(-5)
$\ldots$
\\
2 & 8 & -2.229703572181493732352240313021672944986858675507693036222976(-5)
$\ldots$
\\
\hline
3 & 2 & -1.092764452688696718233957044460372231874277602428901489109438(-1)
$\ldots$
\\
3 & 3 & -7.176813165338438143871896571868137568859537620262992178861887(-3)
$\ldots$
\\
3 & 4 & -6.957183997016348677998754615917611673908374611362524278817441(-4)
$\ldots$
\\
3 & 5 & -7.659277012695409306374743110808079101898869002639397513312088(-5)
$\ldots$
\\
3 & 6 & -8.918921902960271370285096859450098504725826286458231248526033(-6)
$\ldots$
\\
3 & 7 & -1.068734610688718635883673494132633669094248021103907624353475(-6)
$\ldots$
\\
3 & 8 & -1.301295684059645175221229018448850667367013695078510484904330(-7)
$\ldots$
\\
\hline
4 & 2 & -3.603726094351798848506626656181111241130836664796955962932855(-2)
$\ldots$
\\
4 & 3 & -1.174116309572987946977816010618872204640816441822335628644463(-3)
$\ldots$
\\
4 & 4 & -5.722998858912958006017304600902369401289250874897789676420392(-5)
$\ldots$
\\
4 & 5 & -3.165566369796062449665250347230914962435474993776593819875806(-6)
$\ldots$
\\
4 & 6 & -1.848763022618552000611470797190861080535660685413458944944226(-7)
$\ldots$
\\
4 & 7 & -1.109730619419583432307793855949753251494307851245851961690094(-8)
$\ldots$
\\
4 & 8 & -6.763771157059229099811578675678873194164157816965106106489574(-10)
$\ldots$
\\
\hline
5 & 2 & -1.102162098505070183104131920053734921658299916521679210977474(-2)
$\ldots$
\\
5 & 3 & -1.806134929387963117989216971707723297041186986588483316621559(-4)
$\ldots$
\\
5 & 4 & -4.431364427680593899920896902337095946811075675701701222137028(-6)
$\ldots$
\\
5 & 5 & -1.230203263758791942696292884312459227227740434277696211594902(-7)
$\ldots$
\\
5 & 6 & -3.599723532708191784184538837522982283120808246510836532161742(-9)
$\ldots$
\\
5 & 7 & -1.081638288290161011509961768322567652158899297605137495433586(-10)
$\ldots$
\\
5 & 8 & -3.298569076163768263274720121521048493116765743835081292530352(-12)
$\ldots$
\\
\hline
6 & 2 & -3.232720312150523118304098243969541303542841741356062914953186(-3)
$\ldots$
\\
6 & 3 & -2.676915386444316803744871335276940423982308402686133969144141(-5)
$\ldots$
\\
6 & 4 & -3.302884682590606781698851083367477685502328729832308044821381(-7)
$\ldots$
\\
6 & 5 & -4.597234953883457819527050674957592810761129310502451413327254(-9)
$\ldots$
\\
6 & 6 & -6.735520381085029586780892400600876281099743200915951990514023(-11)
$\ldots$
\\
6 & 7 & -1.012733294209721415331392073627541785560477001434410873508589(-12)
$\ldots$
\\
6 & 8 & -1.544937368294918275696159095802943943066554391662366563388540(-14)
$\ldots$
\\
\hline
\end{tabular}
\label{tab.Pkpri}
\end{table}
\subsection{M\"obius (Square-free) Variant}
Reduction of the summation to $k$-almost primes with $k$ distinct
prime factors defines a signed variant of the prime zeta functions:
\begin{defn}
\begin{equation}
P_k^{(\mu)}(s)\equiv \sum_{\substack{n=1\\ \Omega(n)=k}}^\infty \frac{\mu(n)}{n^s}
=
(-1)^k \negthickspace\sum_{\substack{n=1\\ \Omega(n)=\omega(n)=k}}^\infty \frac{1}{n^s}
;\quad
P_1^{(\mu)}(s)=-P(s)
,
\end{equation}
where $\omega(.)$ is the number of distinct prime factors of its argument.
\end{defn}
\begin{rem}
The criterion $\Omega(n)=\omega(n)=k$ selects the
prime number products (square-free $k$-almost primes) of
A000040 ($k=1$),
A006881 ($k=2$), A007304 ($k=3$), A046386 ($k=4$), A046387 ($k=5$),
A067885 ($k=6$), A123321 ($k=7$),
A123322 ($k=8$),
A115343 ($k=9$),
etc.
\end{rem}
The sum
\begin{equation}
\frac{1}{\zeta(s)}=\sum_{n=1}^\infty \frac{\mu(n)}{n^s}=\sum_{k=1}^\infty P_k^{(\mu)}(s)
\end{equation}
converges for $s>\frac{1}{2}$ if the Riemann hypothesis holds \cite{ErdelyiIII}.
Application of the multinomial expansion \cite[(24.1.2)]{AS} to the powers $P^k(s)$ leads to
the recurrences
\begin{thm}
\begin{gather}
P_2^{(\mu)}(s) = 
\frac{P^2(s)-P(2s)}{2!};
\\
-P_3^{(\mu)}(s) = 
\frac{P^3(s)-3P(s)P(2s)+2P(3s)}{3!};
\\
P_4^{(\mu)}(s) = 
\frac{P^4(s)-6P^2(s)P(2s)+3P^2(2s)+8P(s)P(3s)-6P(4s)}{4!};
\end{gather}
\begin{multline}
P_k^{(\mu)}(s)
= \frac{1}{k!}
\sum_{\substack{k_1+2k_2+3k_3+\dotsb+kk_k=k\\ k_k\ge 0}}
(-1)^{m}
(k;k_1k_2\ldots k_k)^*
\\
\times
P^{k_1}(s)P^{k_2}(2s)\dotsm P^{k_k}(ks)
,
\label{eq.Pmobks}
\end{multline}
where
$m\equiv k_1+k_2+k_3+\dotsb+k_k$.
\end{thm}
Redistributing the sign with
\begin{equation}
(-1)^m
P^{k_1}(s)P^{k_2}(2s)\dotsm P^{k_k}(ks)
=
P^{(\mu)k_1}_1(s)P^{(\mu)k_2}_1(2s)\dotsm P^{(\mu)k_k}_1(ks),
\end{equation}
shows that a recurrence equivalent to (\ref{eq.cyclRec}) is applicable.
\begin{rem}
Sums over odd indices are \cite{SebahGourdonPi}:
\begin{gather}
\sum_{k=1}^\infty P^{(\mu)}_{2k-1}(2)=-\frac{9}{2\pi^2};\quad
\sum_{k=1}^\infty P^{(\mu)}_{2k-1}(2s)=-\frac{\zeta^2(2s)-\zeta(4s)}{2\zeta(2s)\zeta(4s)}.
\end{gather}
\end{rem}

\begin{table}
\caption{Almost-prime zeta functions $P_k^{(\mu)}(s)$ at small integer arguments $s$
computed from (\ref{eq.Pmobks}).
}
\begin{tabular}{|r|r|l|}
\hline
$k$ & $s$ & $P_k^{(\mu)}(s)$ \\
\hline
2 & 2 & 6.37672945847765432801316294807193836128782162900370736592109(-2)
$\ldots$
\\
2 & 3 & 6.73594662213544672456228258677680141934623660580421211246428(-3)
$\ldots$
\\
2 & 4 & 9.33269102119509074753434906896208799524913336159126727476407(-4)
$\ldots$
\\
2 & 5 & 1.42408850419374548659705490588033249391875510740905235597858(-4)
$\ldots$
\\
2 & 6 & 2.26806975268639528950901474490066255999464767652487721194642(-5)
$\ldots$
\\
2 & 7 & 3.68874503144697741103258077849260322707848294774819695331204(-6)
$\ldots$
\\
2 & 8 & 6.06493665920244704921794541628232570617319835668401899600473(-7)
$\ldots$
\\
\hline
3 & 2 & -3.6962441634528353783955346323946681155915397130304272497472(-3)
$\ldots$
\\
3 & 3 & -6.6148651246349939521729829639111115641021894727404106069829(-5)
$\ldots$
\\
3 & 4 & -1.7271458093722304630212588271041732572671841808663978769400(-6)
$\ldots$
\\
3 & 5 & -5.0940194598826356852005108113237778072664477868404106866961(-8)
$\ldots$
\\
3 & 6 & -1.5823229154549293389076239250147682789572853572784414723597(-9)
$\ldots$
\\
3 & 7 & -5.0453603114670647581939240532674248065812561949214992614300(-11)
$\ldots$
\\
3 & 8 & -1.6329431236938215954403416483738501457291953254025603265554(-12)
$\ldots$
\\
\hline
4 & 2 & 1.05117508492309807485233009466098526324680558243958672947068(-4)
$\ldots$
\\
4 & 3 & 2.14173193213549705893739943930728906255490278218470044772058(-7)
$\ldots$
\\
4 & 4 & 7.29603168874401925790854604647164607340676506339572130564981(-10)
$\ldots$
\\
4 & 5 & 2.96196721173369084753821609237130107261806114215416577935264(-12)
$\ldots$
\\
4 & 6 & 1.29842711892568424473824206373082938247800573758362639862211(-14)
$\ldots$
\\
4 & 7 & 5.91005736941452577777835874048660952776699696196276661514317(-17)
$\ldots$
\\
4 & 8 & 2.74348813375914336305677937087519697186480294957450635783864(-19)
$\ldots$
\\
\hline
5 & 2 & -1.6620822035796812822471192427038246964759250664122544710020(-6)
$\ldots$
\\
5 & 3 & -2.6408478825460477590567284836043051634661147998573638254764(-10)
$\ldots$
\\
5 & 4 & -7.6296745513152797837319954330917616771026320680788925154580(-14)
$\ldots$
\\
5 & 5 & -2.6674855904890302723828983961667330621935895399254169680514(-17)
$\ldots$
\\
5 & 6 & -1.0124771878309455119843894285163727048400779511505829601925(-20)
$\ldots$
\\
5 & 7 & -4.0089171582464324580016512585964432909883144408807402330430(-24)
$\ldots$
\\
5 & 8 & -1.6267872202198776712860024734866969196014335243663121425423(-27)
$\ldots$
\\
\hline
6 & 2 & 1.61508116616705485066326829316589367443121196123113678467541(-8)
$\ldots$
\\
6 & 3 & 1.53530343588928080283456861104632555466175700586234793665716(-13)
$\ldots$
\\
6 & 4 & 2.99324100743804909283239785893475963049642868205678987129596(-18)
$\ldots$
\\
6 & 5 & 7.41916124138504344746395187399890095182210117294534971585873(-23)
$\ldots$
\\
6 & 6 & 2.05602063767548079577453001329402669576701618478791822725187(-27)
$\ldots$
\\
6 & 7 & 6.06313176438106233958227084211400755691448415132455986317069(-32)
$\ldots$
\\
6 & 8 & 1.85800289034470997723167773802340218486625624481826862206504(-36)
$\ldots$
\\
\hline
\end{tabular}
\label{tab.Pkmu}
\end{table}

\section{Applications} \label{sec.4}
\subsection{Geometric Series}
By integrating \cite[(1.511)]{GR}
\begin{equation}
\log(1-x)= -\sum_{n=1}^\infty \frac{x^n}{n}
\label{eq.logx1}
\end{equation}
we obtain the generating function \cite[(1.513.5)]{GR},
\begin{equation}
x+(1-x)\log(1-x)=\sum_{n=2}^\infty \frac{x^n}{n(n-1)}
.
\label{eq.logx2}
\end{equation}
In the limit $x\to 1$ we find \cite[(0.141)]{GR}
\begin{equation}
\sum_{n=2}^\infty \frac{1}{n(n-1)}=1
.
\label{eq.Nk1}
\end{equation}
Larger powers of the first
factor of the denominator display partial sums of zeta functions through
partial fraction decomposition,
\begin{equation}
\sum_{n=2}^\infty \frac{1}{n^2(n-1)}
=
\sum_{n=2}^\infty \left(\frac{1}{n(n-1)}-\frac{1}{n^2}\right)
=
2-\frac{\pi^2}{6};
\end{equation}
\begin{equation}
\sum_{n=2}^\infty \frac{1}{n^3(n-1)}
=
\sum_{n=2}^\infty \left(\frac{1}{n(n-1)}-\frac{1}{n^2}-\frac{1}{n^3}\right)
=
3-\frac{\pi^2}{6}-\zeta(3);
\end{equation}
\begin{equation}
\sum_{n=2}^\infty \frac{1}{n^s(n-1)}
=
s-
\sum_{l=2}^s
\zeta(l); \quad s\ge 1
.
\label{eq.nmmpf}
\end{equation}
Examples of these numbers are collected in Table \ref{tab.mm},
A152416 and A152419. Their sum is
\begin{equation}
\sum_{s=1}^\infty \sum_{n=2}^\infty \frac{1}{n^s}\,\frac{1}{n-1}
=
\sum_{n=2}^\infty \frac{1}{(n-1)^2}=\zeta(2)=\frac{\pi^2}{6}.
\end{equation}
\begin{table}
\caption{Series of the form (\ref{eq.nmmpf}).
}
\begin{tabular}{|l|l|}
\hline
$s$ & $\sum_{n=2}^\infty 1/[n^s(n-1)]$ \\
\hline
1 & 1.
\\
2 & 3.550659331517735635275848333539748107810500987932015622(-1)
$\ldots$\\3 & 1.530090299921792781278466718425248200160638064527026804(-1)
$\ldots$\\4 & 7.068579628104108661184297530135691724131285453397577278(-2)
$\ldots$\\5 & 3.375804113767116028047748884432274918423193503206296081(-2)
$\ldots$\\6 & 1.641497915322202056595955905340222128241444499920939897(-2)
$\ldots$\\7 & 8.065701771299193726162009203605461682550884433970692555(-3)
$\ldots$\\8 & 3.988345573354854347476770694952996423590093784120672226(-3)
$\ldots$\\9 & 1.979952747272639929624001462540935937984242389231915677(-3)
$\ldots$\\10 & 9.853776194545545924780425622219189319647108247543984198(-4)
$\ldots$
\\
\hline
\end{tabular}
\label{tab.mm}
\end{table}

\begin{defn}
Restriction of the summation in (\ref{eq.nmmpf}) to $k$-almost primes defines
a
set of constants $B_{k,s}$ \cite{GrahMonMM118},
\begin{equation}
B_{k,s}\equiv 
\sum_{\substack{n=2\\ \Omega(n)=k}}^\infty \frac{1}{n^s(n-1)}
=
\sum_{\substack{n=2\\ \Omega(n)=k}}^\infty \sum_{l=0}^\infty \frac{1}{n^{s+1+l}}
=
\sum_{l=0}^\infty
P_k(s+1+l)
.
\label{eq.Nksdef}
\end{equation}
\end{defn}
This reduction to a geometric series and sum over the $P_k$ has been used
to calculate Table \ref{tab.Nk1}.
\begin{table}
\caption{Some values of $B_{k,1}$.
In accordance with (\ref{eq.Nk1}), the series limit of the partial sums is
1 as $k\to \infty$.
}
\begin{tabular}{|l|l|}
\hline
$k$ & $B_{k,1}=\sum_{n,\Omega(n)=k} 1/[n(n-1)]$ \\
\hline
1 &            .77315666904979512786436745985594239561874133608318604831100606
$\ldots$\\2 &  .17105189297999663662220256437237421399124661203550059749107997
$\ldots$\\3 &  .41920339281764199227805032233471158322784525420828606710238790(-1)
$\ldots$\\4 &  .10414202346301156141109353888171559234184072973208943335673068(-1)
$\ldots$\\5 &  .25944317032356863609340108179412019910406474149863463566912649(-2)
$\ldots$\\6 &  .64712678336846601104554817217635310331423959328614964903156640(-3)
$\ldots$\\7 &  .16154547889045106884023528793253084539703632404976961733195128(-3)
$\ldots$\\8 &  .40350403394466614988860237144035458196679194641891917284153345(-4)
$\ldots$\\9 &  .10082343610557897391498490786448232831594447756388459395195738(-4)
$\ldots$\\10 & .25198413274347214707213045392269003455525827323722877742607328(-5)
$\ldots$\\11 & .62985737261498933173999701960384580565617680181013417876197437(-6)
$\ldots$\\12 & .15745036385232517679881727385782020683287536986374137174089714(-6)
$\ldots$\\13 & .39360719611599520681959076312200454332371020210585755644656168(-7)
$\ldots$\\14 & .98399321710455906992308452734477378441132531283484928659273734(-8)
$\ldots$\\15 & .24599505388932024146078978227171922132133360668716211893213680(-8)
$\ldots$\\16 & .61498340075560680170866561648700452467948343944318880698907151(-9)
$\ldots$\\17 & .15374530193448859466794967919224844842656125026606949920935630(-9)
$\ldots$\\18 & .38436254839132442407411486355109175383764384816988280391113599(-10)
$\ldots$\\19 & .96090546444389183823411185748926311892868554398657649851752505(-11)
$\ldots$\\20 & .24022625018522170257126899356289067549236119426157474416944241(-11)
$\ldots$\\21 & .60056547765687746047761256134140717295823654635402409654484759(-12)
$\ldots$\\22 & .15014135061632779005230691663603510295794137329809014913244327(-12)
$\ldots$\\23 & .37535335268553839790798447165432785944777132136305316859771955(-13)
$\ldots$\\24 & .93838335149705572793424612556180553285473359587783011738232439(-14)
$\ldots$\\25 & .23459583405297725332944711004418618278736884444821387508950778(-14)
$\ldots$
\\
\hline
\end{tabular}
\label{tab.Nk1}
\end{table}
This definition introduces an analog to (\ref{eq.zetpart}),
\begin{equation}
\sum_{n=2}^\infty \frac{1}{n^s(n-1)}
=
\sum_{k=1}^\infty B_{k,s}
.
\label{eq.mm}
\end{equation}
\begin{rem}
$B_{1,1}$ is A136141, calculated by Cohen \cite{Cohen}.
$B_{1,2}$ is A152441.
$B_{2,1}$ is A152447.
\end{rem}

Projection of (\ref{eq.nmmpf}) onto the $n$  of a fixed $\Omega(n)$ yields
\begin{equation}
B_{k,s}
=
B_{k,1}
-
\sum_{l=2}^s
P_k(l)
.
\end{equation}
The benefit of this formula is that the $B_{k,s}$ can be derived from $B_{k,1}$
without accumulating the individual terms of the geometric series
proposed in (\ref{eq.Nksdef}), reaching back
to Tables \ref{tab.P} and \ref{tab.Pk} instead.
Consider for example $B_{3,2}
\approx 0.003404
=B_{3,1}-P_3(2)
\approx 0.041920-0.038516$
in Table \ref{tab.Nks}.

\begin{table}
\caption{Some values of $B_{k,s}$.
The series limits of the partial sums $\sum_k B_{k,s}$
are in Table \ref{tab.mm}.
}
\begin{tabular}{|l|l|l|}
\hline
$k$ & $s$ & $B_{k,s}=\sum_{\Omega(n)=k} 1/[n^s(n-1)]$ \\
\hline
1 & 2 &            .32090924900872962935782409502369446144550999284329362657458713
$\ldots$\\2 & 2 &  .30291458630973248399451638958496960216327336030620333566930993(-1)
$\ldots$\\3 & 2 &  .34041462990695583149671096054107628838843579424470347730472318(-2)
$\ldots$\\4 & 2 &  .40476614481790532069851037008630456765933446284259484392905784(-3)
$\ldots$
\\
\hline
          1 & 3 &  .14614660970928609293471078035798776047009787091714433668591512
$\ldots$\\2 & 3 &  .64854251582012887207556831056607595873951575502718767213033361(-2)
$\ldots$\\3 & 3 &  .35478421673524536821901075833145988403566339382745677940994463(-3)
$\ldots$\\4 & 3 &  .20861610202178413235709527570020237570255452209283788001114446(-4)
$\ldots$\\5 & 3 &  .12661340580586229820433777990228912341234438356263326260613819(-5)
$\ldots$
\\
\hline
          1 & 4 &  .69153469945039247992091484424829890308056811202301146420977111(-1)
$\ldots$\\2 & 4 &  .14907506895639490854788090560814702648930997020410089927942399(-2)
$\ldots$\\3 & 4 &  .40356719902303626036886094140740676728492698470022769157503418(-4)
$\ldots$\\4 & 4 &  .11819765739964985563009679542861280192566527050014484191210481(-5)
$\ldots$\\5 & 4 &  .35812310330127538835013908172912057990127682625080625278844028(-7)
$\ldots$
\\
\hline
          1 & 5 &  .33398452461114990859273241885974179176359534475649814730884436(-1)
$\ldots$\\2 & 5 &  .35473826470759431896325286534569759919935164601490043664281503(-3)
$\ldots$\\3 & 5 &  .47794665316209349180059179127476174264255324471931883104976848(-5)
$\ldots$\\4 & 5 &  .69871075602351514235551110353718404446313416171798954710770233(-7)
$\ldots$\\5 & 5 &  .10577117291999709417355375910626719164172948014247187268185467(-8)
$\ldots$
\\
\hline
          1 & 6 &  .16328365610478477905139568619914779966773592601105569997721067(-1)
$\ldots$\\2 & 6 &  .86031097146184686541514125731610045619472901302328472349204898(-4)
$\ldots$\\3 & 6 &  .57819099766026370361808361682441463146581149531325986267386085(-6)
$\ldots$\\4 & 6 &  .42224059396278682956248009245427510821206201028665096074007060(-8)
$\ldots$\\5 & 6 &  .31946136165040413132735812808230110887705358254036598872046674(-10)
\\
\hline
\end{tabular}
\label{tab.Nks}
\end{table}

The square-free variant of (\ref{eq.Nksdef}) is
\begin{defn}
\begin{equation}
B_{k,s}^{(\mu)}\equiv 
\sum_{\substack{n=2\\ \Omega(n)=k}}^\infty \frac{\mu(n)}{n^s(n-1)}
.
\label{eq.Bksmudef}
\end{equation}
\end{defn}
As in (\ref{eq.Nksdef}), there is a representation generated by expansion as a geometric
series, and another one from the decomposition in partial fractions:
\begin{equation}
B_{k,s}^{(\mu)}
=
\sum_{l=0}^\infty
P_k^{(\mu)}(s+1+l)
=
B_{k,1}^{(\mu)}
-
\sum_{l=2}^s
P_k^{(\mu)}(l)
.
\end{equation}
Explicit values follow in Tables \ref{tab.Nk1mu} and \ref{tab.Nksmu}. The special
values of
\begin{equation}
B_{1,s}^{(\mu)}= -B_{1,s}
\end{equation}
can be read off Table \ref{tab.Nk1} and \ref{tab.Nks}.
Because the squared primes are those $2$-almost primes which are not square-free,
constants like
\begin{equation}
\sum_{p} \frac{1}{p^{2s}(p^2-1)}
=
\sum_{l=0}^\infty P(2(1+s+l))
=B_{2,s}-B_{2,s}^{(\mu)}
\end{equation}
can be extracted subtracting values of Tables \ref{tab.Nks} and \ref{tab.Nksmu}.
\begin{table}
\caption{Some values of (\ref{eq.Bksmudef}) at $s=1$.
}
\begin{tabular}{|l|l|}
\hline
$k$ & $(-1)^kB_{k,1}^{(\mu)} = \sum_{\Omega(n)=\omega(n)=k} 1/[n(n-1)]$ \\
\hline
1 &             .77315666904979512786436745985594239561874133608318604831100606
$\ldots$\\2 &   .71606015364062950689014905233278570032977577496764766996881566(-1)
$\ldots$\\3 &   .37641725351677739987897144642934934884513171678284733095361749(-2)
$\ldots$\\4 &   .10533241426370309073189561057615806202677017536174797996070137(-4)
$\ldots$\\5 &   .16623463646913663848631359812999970758043030193325500480455497(-5)
$\ldots$\\6 &   .16150965195007452635819510602638099761518437046602308615745291(-7)
$\ldots$\\7 &   .10379657945831823210405127184359674418404044704893780297641014(-9)
$\ldots$\\8 &   .46608164350339032665792725000352856059530910663817211356715781(-12)
$\ldots$\\9 &   .15257916508734074179181916995271562612876849081982812066031366(-14)
$\ldots$\\10&   .37674605405462816954865036936087328930035659282686679768832686(-17)
$\ldots$\\11&   .72146307104358813058965067637397207142589348895613246863457289(-20)
$\ldots$\\12&   .10966184934068789212128440266173951091947363648217302666195170(-22)
$\ldots$\\13&   .13491272590180303187806861812072181049977941362569266836182658(-25)
$\ldots$\\14&   .13659203913426921565993991785542457997100724146742659833750046(-28)
$\ldots$\\15&   .11544550299305819020362074848337321347361648310620287413349137(-31)
$\ldots$\\16&   .82468163663946083189906712482357029072394743693281076128206927(-35)
$\ldots$\\17&   .50333187700965107357204271671076913932386822516030741929080100(-38)
$\ldots$\\18&   .26498528956530623344847668001272588778001122254618728853331610(-41)
$\ldots$\\19&   .12135579687796160737893591360857756602310726382140060310426006(-44)
$\ldots$\\20&   .48713430940243727394794760257290243400665149725603953693181101(-48)
$\ldots$\\21&   .17255587625660721457243893666076741199905259181309475423500537(-51)
$\ldots$\\22&   .54270137257558710718858909894832920484581333704208745878297310(-55)
$\ldots$\\23&   .15238881562580861279422356191754120203296455021672076901780409(-58)
$\ldots$\\24&   .38397551449809407896292097416496474894690864603905322011328373(-62)
$\ldots$\\25&   .87221420134797085669777164633365058358980257342379119911641791(-66)
$\ldots$
\\
\hline
\end{tabular}
\label{tab.Nk1mu}
\end{table}
\begin{table}
\caption{Some values of $(-1)^kB_{k,s}^{(\mu)}$.
}
\begin{tabular}{|l|l|l|}
\hline
$k$ & $s$ & $(-1)^kB_{k,s}^{(\mu)}$ \\
\hline
          2 & 2 &   .78387207792864074088832757525591864200993612067276933376705990(-2)
$\ldots$\\3 & 2 &   .67928371714938620394179831898825372859777454798046059788912729(-4)
$\ldots$\\4 & 2 &   .21490577139328324666260111005953570208961711778930701363329944(-6)
$\ldots$\\5 & 2 &   .26416111168510261601673859617237932837795292029557704345022153(-9)
$\ldots$
\\
\hline
          2 & 3 &   .11027741571509606843209931657823850007531246009234812252063101(-2)
$\ldots$\\3 & 3 &   .17797204685886808724500022597142572187555600706419537190829226(-5)
$\ldots$\\4 & 3 &   .73257817973354076886116612880679583412683957083696886124052352(-9)
$\ldots$\\5 & 3 &   .76323430497840111065747811948812031341440309840660902575023325(-13)
$\ldots$
\\
\hline
          2 & 4 &   .16950505503145160956755825888617620122821126476435449772990259(-3)
$\ldots$\\3 & 4 &   .52574659216450409428743432610083961488375889775555842142889903(-7)
$\ldots$\\4 & 4 &   .29750108591388430703115241596312267861630644973967306755420261(-11)
$\ldots$\\5 & 4 &   .26684984687313228427857617894414570413989159871977420443245736(-16)
$\ldots$
\\
\hline
          2 & 5 &   .27096204612077060907852768298142951836335754023449262132044308(-4)
$\ldots$\\3 & 5 &   .16344646176240525767383244968461834157114119071517352759281108(-8)
$\ldots$\\4 & 5 &   .13043647405152222773308067259925713545003355242564896189378205(-13)
$\ldots$\\5 & 5 &   .10128782422925704028633932747239792053264472723250762731720244(-19)
$\ldots$
\\
\hline
          2 & 6 &   .44155070852131080127626208491363262363892772582004900125801062(-5)
$\ldots$\\3 & 6 &   .52141702169123237830700571831415136754126549873293803568336271(-10)
$\ldots$\\4 & 6 &   .59376215895380325925646622617419720223297866728632203157065385(-16)
$\ldots$\\5 & 6 &   .40105446162489087900384620760650048636932117449331297950429083(-23)
\\
\hline
\end{tabular}
\label{tab.Nksmu}
\end{table}

\subsection{Hurwitz Zeta Decompositions}
\begin{defn} The Hurwitz Zeta Function is
\begin{equation}
\zeta(s,a)\equiv \sum_{n=0}^\infty \frac{1}{(n+a)^s}
=
\sum_{n=1}^\infty \frac{1}{(n+a-1)^s}
;\quad \Re s>1,\quad \Re a>1
.
\label{eq.Hzeta}
\end{equation}
\end{defn}
On the trot, we project this sum onto the subspaces of $k$-almost primes, too,
\begin{defn} (Hurwitz Prime and Almost-Prime Zeta Functions)
\begin{equation}
P_k(s,a)\equiv \sum_{\substack{n=2\\ \Omega(n)=k}}^\infty \frac{1}{(a-1+n)^s}
;
\quad
P_k(s,1)=P_k(s).
\label{eq.Pksadef}
\end{equation}
\begin{equation}
P_k^{(\mu)}(s,a)\equiv \sum_{\substack{n=2\\ \Omega(n)=k}}^\infty \frac{\mu(n)}{(a-1+n)^s}
;
\quad
P_k^{(\mu)}(s,1)=P_k^{(\mu)}(s).
\end{equation}
\end{defn}
This partitions (\ref{eq.Hzeta}) into
\begin{equation}
\zeta(s,a)=\frac{1}{a^s}+\sum_{k=1}^\infty P_k(s,a),
\end{equation}
generalizing (\ref{eq.zetpart}).
\begin{rem}
By series expansion \cite[(1.112.1)]{GR}, the sum rule associated with (\ref{eq.Nksdef}) is
\begin{equation}
B_{k,s}
=
\sum_{\substack{n=2\\ \Omega(n)=k}}^\infty \frac{1}{(n-1)^{1+s}}\,\frac{1}{1+\frac{1}{n-1}}
=
\sum_{l=0}^\infty
(-1)^lP_k(1+s+l,0).
\end{equation}
\end{rem}
The reduction of (\ref{eq.Pksadef}) to the Prime Zeta Functions is obtained by the binomial
expansion \cite[(1.110)]{GR}\cite[Entry 22]{BerndtPIACS92}
\begin{multline}
P_k(s,a)
=
\sum_{\substack{n=2\\ \Omega(n)=k}}^\infty \frac{1}{n^s}\,
\frac{1}{
\left(1-\frac{1-a}{n}\right)^s
}
=
\sum_{\substack{n=2\\ \Omega(n)=k}}^\infty \frac{1}{n^s}\,
\sum_{l=0}^\infty
\binom{-s}{l}(-1)^l\frac{(1-a)^l}{n^l}
\\
=
\sum_{\substack{n=2\\ \Omega(n)=k}}^\infty \frac{1}{n^s}\,
\sum_{l=0}^\infty
\frac{(s)_l}{l!}
\frac{(1-a)^l}{n^l}
=
\sum_{l=0}^\infty
\frac{(s)_l}{l!}
(1-a)^l P_k(s+l)
.
\label{eq.PksaofPk}
\end{multline}
\begin{equation}
P_k^{(\mu)}(s,a)
=
\sum_{l=0}^\infty
\frac{(s)_l}{l!}
(1-a)^l P_k^{(\mu)}(s+l)
.
\end{equation}
Pochhammer's symbol $(s)_l\equiv \Gamma(s+l)/\Gamma(s)$ is introduced
to simplify the notation \cite[(6.1.22)]{AS}.
Brute-force
accumulation
of partial sums over $l$ for the case $a=0$ yields Table \ref{tab.H0}.
Resummation in (\ref{eq.Pksadef}) imposes the sum rule
\begin{equation}
\sum_{k=1}^\infty P_k(s,0) = \zeta(s)
.
\end{equation}

\begin{table}
\caption{Some values of $P_{k}(s,0)$.
Cohen has reported $P_1(2,0)$
\cite{Cohen},
which is A086242 \cite{EIS}.
}
\begin{tabular}{|l|l|l|}
\hline
$k$ & $s$ & $P_k(s,0)=\sum_{\Omega(n)=k} 1/(n-1)^s$ \\
\hline
          1 & 2 &  1.37506499474863528791725313052243969917959996017531745870918933
$\ldots$\\
          2 & 2 &  .209788323940019492755368602469189236268613932921851343752817089
$\ldots$\\
          3 & 2 &  .457250649473356179509462896789867962239663545084304304104898823(-1)
$\ldots$\\
          4 & 2 &  .108410864467466394130511466215222764287414954775183405819497437(-1)
$\ldots$\\
          5 & 2 &  .264509027543661792197346769155986673457680151009439096206096105(-2)
$\ldots$\\
\hline
          1 & 3 &  1.14752909775858004693328380628213040164476473552511225527582412
$\ldots$\\
          2 & 3 &  .497610511326665981875950866377952405841340774755637814783395177(-1)
$\ldots$\\
          3 & 3 &  .425728972912996225107947200519463804614525467386937187206041965(-2)
$\ldots$\\
          4 & 3 &  .450337561601661496970504401383245773976032018827147548523156327(-3)
$\ldots$\\
          5 & 3 &  .519996368729267202444682115275828429696087549531802544406653646(-4)
$\ldots$\\
\hline
          1 & 4 &  1.06736011227157169811527402065258703893525859304550836811508241
$\ldots$\\
          2 & 4 &  .144251867050125347206321771385325774655732916871746611337540069(-1)
$\ldots$\\
          3 & 4 &  .511594573334901127268316650383244703997517641565236190813247750(-3)
$\ldots$\\
          4 & 4 &  .248732308184070333068620411510121468878257641611466236659372052(-4)
$\ldots$\\
          5 & 4 &  .138025203867843406309850933144822342517665004454740887033517806(-5)
$\ldots$\\
\hline
          1 & 5 &  1.03237100597834196585177592063868294503482496931776985434695475
$\ldots$\\
          2 & 5 &  .448814553317860895890452521422141457224197062051853943952316742(-2)
$\ldots$\\
          3 & 5 &  .670525710706644570099917906002173405129522585443401784675937527(-4)
$\ldots$\\
          4 & 5 &  .150947118336184838341631597141043715807449152523875014126944996(-5)
$\ldots$\\
          5 & 5 &  .403832718324949534495289315322822843326177902298377065444183828(-7)
$\ldots$\\
\hline
          1 & 6 &  1.01589201139972411006675918457325510103473578053731168636932226
$\ldots$\\
          2 & 6 &  .144181860839192758202651548654521805154418089925454726071029138(-2)
$\ldots$\\
          3 & 6 &  .913523107817850408748985104163510342741650140340861620051453925(-5)
$\ldots$\\
          4 & 6 &  .954942045591842308217096235482438389525995484180595401298115858(-7)
$\ldots$\\
          5 & 6 &  .123322395562635072773966382461892764518565670709235951708582994(-8)
$\ldots$\\
\hline
\end{tabular}
\label{tab.H0}
\end{table}

The derivative of (\ref{eq.PksaofPk}) with respect to $s$ is
\begin{gather}
\frac{\partial}{\partial s}P_k(s,a) = P_k'(s,a) =
- \sum_{\substack{n=2\\ \Omega(n)=k}}^\infty \frac{\log(a-1+n)}{(a-1+n)^s}
\nonumber \\
=
\sum_{l=0}^\infty
\frac{(s)_l(1-a)^l}{l!}
\left[
(\psi(s+l)-\psi(s)) P_k(s+l) + P_k'(s+l)
\right]
\label{eq.Hlog}
\end{gather}
in terms of digamma functions $\psi$ \cite[(6.3)]{AS}.
Cases with $a=0$ are illustrated by Table \ref{tab.Pprime}.
\begin{table}
\caption{Some absolute values of (\ref{eq.Hlog}) at $a=0$.
}
\begin{tabular}{|l|l|l|}
\hline
$k$ & $s$ & $|P_k'(s,0)|=\sum_{\Omega(n)=k} \log(n-1)/(n-1)^s$ \\
\hline
          1 & 2 &  .412038626948453592989536727886919593108693955993272284789334253
$\ldots$\\
          2 & 2 &  .349406402843729094021858840552328411755268945853567444645947573
$\ldots$\\
          3 & 2 &  .122053877557071229252860135234470821144298661274431502394539860
$\ldots$\\
          4 & 2 &  .381324604512840984602505783616974026559738868834901273236913385(-1)
$\ldots$\\
          5 & 2 &  .113453301459081906536851293804267570546096249239456318213977158(-1)
$\ldots$\\
\hline
          1 & 3 &  .122491994469611894418110664126546306148971700364574609076227728
$\ldots$\\
          2 & 3 &  .647147749533078072180333326031334769093826259630394154513521204(-1)
$\ldots$\\
          3 & 3 &  .935118212996044867853444829503103960682129910900028893869639531(-2)
$\ldots$\\
          4 & 3 &  .134268150501925211286561753670177005489358798054445521080127694(-2)
$\ldots$\\
          5 & 3 &  .193322938011967430559782708011999395264147756118102143367684017(-3)
$\ldots$\\
\hline
          1 & 4 &  .505703011282046113580449366113852907995938099426370524892239106(-1)
$\ldots$\\
          2 & 4 &  .172097052934328020906454081592687823346699820793905416728793462(-1)
$\ldots$\\
          3 & 4 &  .105543731770635366656467702388193888381593405142250008490173220(-2)
$\ldots$\\
          4 & 4 &  .705255178358882713248498070972826082493637781890072290188334302(-4)
$\ldots$\\
          5 & 4 &  .492194845815858173046505381936518534428557359827933872744842926(-5)
$\ldots$\\
\hline
          1 & 5 &  .232833728972359609747954746653486493466064393424300179653066541(-1)
$\ldots$\\
          2 & 5 &  .515141863816847304246565784890839865650309296489970884531366999(-2)
$\ldots$\\
          3 & 5 &  .134651494862962805588752454576895989572879746878237372363930174(-3)
$\ldots$\\
          4 & 5 &  .419081737714005222851533545272768768198324374375815699071937904(-5)
$\ldots$\\
          5 & 5 &  .141550001290993202090903719190199844283347760149056451512440612(-6)
$\ldots$\\
\hline
          1 & 6 &  .112107071669452487391848834784015585352560002552358622406258928(-1)
$\ldots$\\
          2 & 6 &  .162309029366737494046263910398940347901702033172772596446220636(-2)
$\ldots$\\
          3 & 6 &  .181009248757925488622702113921921885223627925005301730877648970(-4)
$\ldots$\\
          4 & 6 &  .262385322411092391533263049109867509130968921137179066906794859(-6)
$\ldots$\\
          5 & 6 &  .428617889403977448189220166642396011211315801456016154892551835(-8)
$\ldots$\\
\hline
\end{tabular}
\label{tab.Pprime}
\end{table}

\subsection{Logarithmic Functions}
Repeated integration of (\ref{eq.logx2}) with respect to $x$ yields
a ladder of functions with increasing order $l+1$ of
the polynomial of $n$ in the denominator \cite[(1.513.7)]{GR}\cite{LehmerAMM92},
\begin{gather}
x+(1-x)^2\ln(1-x)
=
\frac{3x^2}{2}
-2
\sum_{n=1}^\infty \frac{x^{n+2}}{n(n+1)(n+2)};
\label{eq.logx3}
\\
x+(1-x)^3\ln(1-x)
=
\frac{5x^2}{2}
-\frac{11x^3}{6}
+6
\sum_{n=1}^\infty \frac{x^{n+3}}{n(n+1)(n+2)(n+3)};
\\
x+(1-x)^4\ln(1-x)
=
\frac{7x^2}{2}
-\frac{13x^3}{3}
+\frac{25x^4}{12}
-24
\sum_{n=1}^\infty \frac{x^{n+4}}{n(n+1)(n+2)(n+3)(n+4)};
\\
x+(1-x)^l\ln(1-x)
=
\sum_{i=2}^l\tau_{i,l}x^i-(-1)^l l!\sum_{n=1}^\infty \frac{x^{n+l}}{n(n+1)\cdots (n+l)}.
\label{eq.logx4}
\end{gather}
The rational coefficients of the auxiliary polynomials obey the recurrence
\begin{equation}
\tau_{2,l}=l-\frac{1}{2},\quad l\ge 2;\quad i\tau_{i,l}= (-1)^i\binom{l-1}{i-1}-l\tau_{i-1,l-1},\quad i\ge 3.
\end{equation}
$i!\tau_{i,l}$ are
tabulated in A067176 and A093905.

Substitution of inverse $k$-almost primes for $x$ in (\ref{eq.logx4}) followed
by summation
defines
constants $L_{k,l}$ which
generalize (\ref{eq.L11}),
\begin{gather}
L_{k,l}\equiv 
\sum_{\Omega(n)=k}\left[\frac{1}{n}+\left(1-\frac{1}{n}\right)^l\ln\left(1-\frac{1}{n}\right)\right]
\nonumber
\\
=
\sum_{s=2}^l\tau_{s,l}P_k(s)
-(-1)^l l!\sum_{s=1}^\infty \frac{1}{s(s+1)\cdots (s+l)}P_k(s+l),
\label{eq.LklDef}
\end{gather}
placed in Table \ref{tab.Lkl}.
At $l=1$, partial fraction decomposition of $1/[s(s+1)]$ in (\ref{eq.LklDef})
implies that $L_{1,1}$ plus the value at $u=1$ in Table \ref{tab.Pmom} equals
the value in (\ref{eq.pmoms}).
Series of the form $\sum_{s\ge 1} P_k(s+l)/\prod_{i=0}^l (s+i)$
are quickly accessible by subtracting the $\tau$-dependent terms
of (\ref{eq.LklDef}) with the aid of Table \ref{tab.Pk}.
Varieties with non-contiguous factors in the denominator
can be reduced to the $L_{k,l}$ basis
by partial fraction synthesis: examples of this approach with one missing
factor $s+2$ or two missing factors $s$ and $s+2$ on the left hand side are
\begin{equation}
\frac{1}{s(s+1)(s+3)}=\frac{1}{s(s+1)(s+2)}-\frac{1}{s(s+1)(s+2)(s+3)}
\end{equation}
and
\begin{equation}
\frac{1}{(s+1)(s+3)}=
\frac{1}{s(s+1)}-
\frac{3}{s(s+1)(s+2)}+\frac{3}{s(s+1)(s+2)(s+3)}.
\end{equation}

\begin{table}
\caption{Logarithmic sums of $k$-almost primes defined in (\ref{eq.LklDef}).
}
\begin{tabular}{|l|l|l|}
\hline
$k$ & $l$ & $L_{k,l}$ \\
\hline
1 &  1& 0.26434004176002164673174824590306461492528050270047150992644
$\ldots$
\\
1 &  2& 0.61210909600861648810973500903074917744536817685379132220140
$\ldots$
\\
1 &  3& 0.83162031622701364585608159976581294069621021299936900413624
$\ldots$
\\
1 &  4& 0.97814589475979949499502571065792051537852172306539544173835
$\ldots$
\\
\hline
2 &  1& 0.07483176753552750888163142884027057172300662738030777281173
$\ldots$
\\
2 &  2& 0.20274600679905285623473576314229971108548607276186973818502
$\ldots$
\\
2 &  3& 0.30956716856446670438488608714185233323846931207915051554876
$\ldots$
\\
2 &  4& 0.39967077199200952417631057083019796152634867816759135315832
$\ldots$
\\
\hline
3 &  1& 0.01979445800870394214182491696866714839521363828266446266352
$\ldots$
\\
3 &  2& 0.05673037199050316150740109374211659920215472144741775457277
$\ldots$
\\
3 &  3& 0.09078044492621127414186177360393579101033447804508695808347
$\ldots$
\\
3 &  4& 0.12224056980762694782137972349781919346596885467095976821370
$\ldots$
\\
\hline
\end{tabular}
\label{tab.Lkl}
\end{table}

\begin{rem}

The case of $x=1/2$ in (\ref{eq.logx1}) is decomposable
in $k$-almost primes as well \cite{SebahGourdonLn2,ChamberlJIS6}
\begin{equation}
\log 2
=
\sum_{n=1}^\infty \frac{1}{n2^n}
= \frac{1}{2}+
\sum_{k=1}^\infty \sum_{\substack{n=1\\ \Omega(n)=k}}^\infty \frac{1}{n2^n}
.
\end{equation}
The dominating components are
\begin{gather}
\sum_{\substack{n=1\\ \Omega(n)=k}}^\infty
\frac{1}{n2^n}
\\
=
0.174087071760979362471993316621554442658749500081033068401\ldots
(k=1)
\nonumber
\\
=
 0.185502662799497065892654852882047774301689318692751270328(-1)\ldots
(k=2)
\nonumber
\\
=
 0.508886353469978309538645063262227959518399790876499567910(-3)\ldots
(k=3)
\nonumber
\\
=
0.956158270230416532955864299570858933417801804520522813509(-6)\ldots
(k=4)
\nonumber
\end{gather}
These sums generated from $|x|< 1$ converge rapidly on their own and do
not need the acceleration strategies outlined in the first sections.
\end{rem}

\begin{rem}
The limits $x\to 1$ in (\ref{eq.logx3})--(\ref{eq.logx4}) generalize
(\ref{eq.Nk1})
\cite[(p.\ 42)]{ZwillingerCRC}
\begin{eqnarray}
\frac{1}{k\cdot k!}
&=&
\sum_{n=1}^\infty\frac{1}{n(n+1)(n+2)\cdots(n+k)}.
\end{eqnarray}
Restricted summation over $k$-almost primes
in the spirit of (\ref{eq.Nksdef})
would split these into another basic type of constants.
\end{rem}

\section{Summary}
Almost-prime zeta functions have been defined by restriction of the
summation of the standard definition of zeta functions to $k$-almost primes.
Their values can be bootstrapped from a multinomial overlay of the values of the 
ordinary prime zeta functions. Efficient schemes to
compute the latter have been employed to calculate series
summed over $k$-almost primes
of some basic inverse polynomials
to high accuracy.

\bibliographystyle{amsplain}
\bibliography{all}

\end{document}